\documentclass[letterpaper, 10 pt, conference]{ieeeconf}  

\IEEEoverridecommandlockouts                              

\title{\LARGE \bf
Rigid Body Motion Estimation based on the Lagrange-d'Alembert Principle
}

\author{Maziar Izadi$^{1}$, Amit K. Sanyal$^{2,\dag}$, Ernest Barany$^{3}$ and Sasi P. Viswanathan$^{1}$
\thanks{$^{1}$M. Izadi and S.P. Viswanathan are with the Department of Mechanical and Aerospace Engineering, New Mexico State University, Las Cruces, NM 88003 USA.
        {\tt\small \{mi,sashi\}@nmsu.edu}}%
\thanks{$^{2}$A. K. Sanyal is with the Department of Mechanical and Aerospace Engineering, Syracuse University, Syracuse, NY 13244 USA.
        {\tt\small aksanyal@syr.edu}}%
\thanks{$^{3}$E. Barany is with the Department of Mathematical Sciences, New Mexico State University, Las Cruces, NM 88003 USA.
        {\tt\small ebarany@nmsu.edu}}%
\thanks{$^{\dag}$Address all correspondence to this author.}
}

\usepackage{graphicx}          
\usepackage{amsmath,amsfonts,amssymb}
\usepackage{mathtools, textcomp}
\usepackage{mathrsfs}
\usepackage{graphicx}
\usepackage[dvips]{epsfig}
\usepackage{epstopdf}
\usepackage{color}
\usepackage{dsfont}

\newcommand{\SO}{\ensuremath{\mathrm{SO(3)}}}

\newcommand{\Ta}{\ensuremath{\mathrm{T}}}
\newcommand{\T}{^{\mbox{\small T}}}
\newcommand{\so}{\ensuremath{\mathfrak{so}(3)}}
\newcommand{\SE}{\ensuremath{\mathrm{SE(3)}}}

\newcommand{\bR}{\ensuremath{\mathbb{R}}}
\newcommand{\bS}{\ensuremath{\mathbb{S}}}

\newcommand{\mrm}{\mathrm}

\newcommand{\diag}{\mbox{diag}}
\newcommand{\bbm}{\begin{bmatrix}}
\newcommand{\ebm}{\end{bmatrix}}
\newcommand{\matl}{\left[ \begin{array}}
\newcommand{\matr}{\end{array} \right]}
\newcommand{\be}{\begin{equation}}
\newcommand{\ee}{\end{equation}}
\newcommand{\bea}{\begin{eqnarray}}
\newcommand{\eea}{\end{eqnarray}}
\newcommand{\beas}{\begin{eqnarray*}}
\newcommand{\eeas}{\end{eqnarray*}}
\newcommand{\nn}{\nonumber}
\newcommand{\mC}{\mathcal{C}}

\newcommand{\cL}{\mathcal{L}}
\newcommand{\cI}{\mathcal{I}}

\newcommand{\cT}{\mathcal{T}}
\newcommand{\cU}{\mathcal{U}}

\newcommand{\di}{\mathrm{d}}

\newcommand{\tr}{\mathrm{trace}}

\newcommand{\lan}{\langle}
\newcommand{\ran}{\rangle}

\newcommand{\cS}{\mathcal{S}}

\newcommand{\ad}[1]{{\mathrm{ad}_{#1}}}          			
\newcommand{\adast}[1]{{\mathrm{ad}_{#1}^\ast}}  			
\newcommand{\Ad}[1]{{\mathrm{Ad}_{#1}}}  			
\DeclareMathOperator{\expm}{{expm}}  			

\newtheorem{theorem}{Theorem}[section]

\def\thesection{\arabic{section}}

\newcommand{\bi}{\begin{itemize}}
\newcommand{\ei}{\end{itemize}}

\DeclareMathAlphabet{\mathpzc}{OT1}{pzc}{m}{it}

\newcommand{\bJ}{\mathbb{J}}

\newcommand{\mpz}{\mathpzc}
\newcommand{\msg}{\mathsf{g}}
\newcommand{\msh}{\mathsf{h}}

\newcommand{\sS}{\mathsf{S}}
\newcommand{\sO}{\mathsf{O}}
\newcommand{\bD}{\mathbb{D}}
\begin{document}

\maketitle
\thispagestyle{empty}
\pagestyle{empty}

\begin{abstract}
Stable estimation of rigid body pose and velocities from noisy measurements, without any 
knowledge of the dynamics model, is treated using the Lagrange-d'Alembert principle from 
variational mechanics.
With body-fixed optical and inertial sensor measurements, a Lagrangian is obtained as the 
difference between a kinetic energy-like term that is quadratic in velocity estimation error and 
the sum of two artificial potential functions; one obtained from a generalization of Wahba's 
function for attitude estimation and another which is quadratic in the position estimate error. An 
additional dissipation term that is linear in the velocity estimation error is introduced, and the 
Lagrange-d'Alembert principle is applied to the Lagrangian with this dissipation. This estimation scheme is discretized using discrete variational mechanics. The presented pose estimator requires optical measurements 
of at least three inertially fixed landmarks or beacons in order to estimate instantaneous pose. 
The discrete estimation scheme can also estimate velocities from such optical measurements. In the presence of bounded measurement noise in the vector measurements, 
numerical simulations show that the estimated states converge to a bounded neighborhood 
of the actual states.
\end{abstract}

\section{INTRODUCTION}\label{Sec1}
Estimation of rigid body translational and rotational motion is indispensable for operations of 
spacecraft, unmanned aerial and underwater vehicles. Autonomous state estimation of a rigid 
body based on inertial vector measurement and visual feedback from stationary landmarks, in 
the absence of a dynamics model for the rigid body, is analyzed here. This estimation scheme can enhance the autonomy and reliability of 
unmanned vehicles in uncertain GPS-denied environments. Salient features of this estimation 
scheme are: (1) use of onboard optical and inertial sensors, with or without rate gyros, for 
autonomous navigation; (2) robustness to uncertainties and lack of knowledge of dynamics; 
(3) low computational complexity for easy implementation with onboard processors; (4) proven 
stability with large domain of attraction for state estimation errors;  
and (5) versatile enough to estimate motion with respect to stationary as well as moving objects. 
Robust state estimation of rigid bodies in the absence of complete knowledge of their dynamics, 
is required for their safe, reliable, and autonomous operations in poorly known conditions. In 
practice, the dynamics of a vehicle may not be perfectly known, especially when the vehicle is 
under the action of poorly known forces and moments. The scheme proposed here has a single, 
stable algorithm for the coupled translational and rotational motion of rigid bodies using 
onboard optical (which may include infra-red) and inertial sensors. This avoids the need for 
measurements from external sources, like GPS, which may not be available in indoor, underwater 
or cluttered environments \cite{leishman2014relative}.
 
Attitude estimators using unit quaternions for attitude representation may be {\em unstable 
in the sense of Lyapunov}, unless they identify antipodal quaternions with a single attitude. 
This is also the case for attitude control schemes based on continuous feedback of unit 
quaternions, as shown in~\cite{Bayadi2014almost}. One adverse 
consequence of these unstable estimation and control schemes is that they end up taking 
longer to converge compared with stable schemes under similar initial conditions and initial 
transient behavior. Continuous-time attitude observers and filtering schemes on 
$\SO$ and $\SE$ have been reported in, e.g., \cite{bonmaro09,Khosravian2015observers,mahapf08,rehbinder2003pose,Vas1}. 
These estimators do not suffer from kinematic singularities like estimators using coordinate 
descriptions of attitude, and they do not suffer from unwinding as they do not use 
unit quaternions. The maximum-likelihood (minimum energy) filtering method of 
Mortensen~\cite{Mortensen} was recently applied to attitude estimation, resulting in a nonlinear 
attitude estimation scheme that seeks to minimize the stored ``energy" in measurement 
errors~\cite{ZamPhD}. 
This scheme is obtained by applying Hamilton-Jacobi-Bellman (HJB) theory to 
the state space of attitude motion. Since the HJB equation can only be 
approximately solved with increasingly unwieldy expressions for higher order 
approximations, the resulting filter is only ``near optimal" up to second order. Unlike filtering 
schemes that are based on approximate or ``near optimal" solutions of the HJB equation and 
do not have provable stability, the estimation scheme obtained here is shown to be almost globally asymptotically stable. Moreover, unlike filters 
based on Kalman filtering, the estimator proposed here does not presume any knowledge 
of the statistics of the initial state estimate or the sensor noise. Indeed, for 
vector measurements using optical sensors with limited field-of-view, the probability 
distribution of measurement noise needs to have compact support, unlike standard Gaussian 
noise processes that are commonly used to describe such noisy measurements.

The variational attitude estimator recently appeared in \cite{Automatica}, where it was 
shown to be almost globally asymptotically stable. Some of the advantages of this scheme over 
some commonly used competing schemes are reported in \cite{ICRA2015}. This paper is the   
variational estimation framework to coupled rotational (attitude) and translational motion, as 
exhibited by maneuvering vehicles like UAVs. In such applications, designing separate state 
estimators for the translational and rotational motions may not be effective and may lead to poor 
navigation. Moreover, like 
other vision-inertial navigation schemes~\cite{shen2013vision}, the estimation 
scheme proposed here does not rely on GPS. However, unlike many other vision-inertial 
estimation schemes, the estimation scheme proposed here can
be implemented without any direct velocity measurements. Since rate gyros are usually corrupted by high noise content and bias \cite{Good2013ECC,Good2014ASME,Grip2015}, 
such a velocity measurement-free scheme can result in fault tolerance in the case of faults with rate 
gyros. Additionally, this estimation scheme can be extended to relative pose estimation between 
vehicles from optical measurements, without direct communications or measurements of 
relative velocities.

\section{NAVIGATION USING OPTICAL AND INERTIAL SENSORS} \label{Sec2}
Consider a vehicle in spatial (rotational and translational) motion. 
\begin{figure}[htb!]
	\centering
	\includegraphics[width=0.475\textwidth]{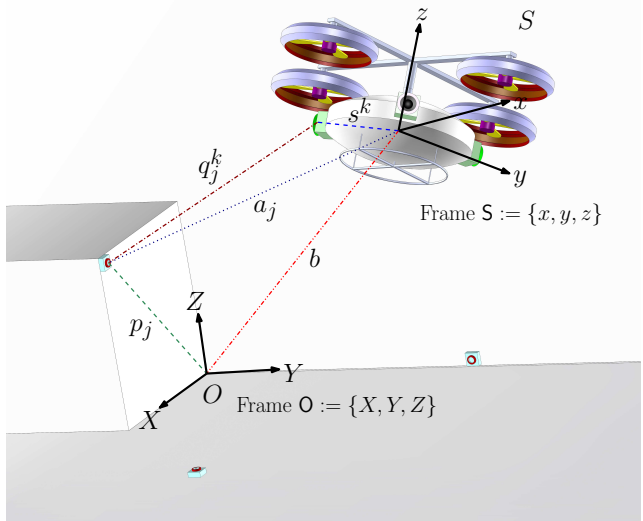}
    	\caption{Inertial landmarks on $O$ as observed from vehicle $S$ with optical  
	measurements.}
    	\label{Frames}
\end{figure}
Onboard estimation of the pose of the vehicle involves assigning a coordinate frame fixed to 
the vehicle body, and another coordinate frame fixed in the environment which takes the role 
of the inertial frame. Let $O$ denote the observed environment and $S$ denote the vehicle. 
Let $\mathsf{S}$ denote a coordinate frame fixed to $S$ and $\sO$ be a coordinate frame 
fixed to $O$, as shown in Fig. \ref{Frames}. Let $R\in\SO$ denote the rotation matrix from frame 
$\sS$ to frame $\sO$ and $b$ denote the position of origin of $\sS$ expressed in frame $\sO$. The pose (transformation) from body fixed frame $\sS$ to inertial frame $\sO$ is then given by
\begin{align} 
&\msg=\bbm R \;\;\;& b\\ 0 \;\;\;& 1\ebm\in\SE.
\label{gDef} 
\end{align}
Consider vectors known in inertial frame $\sO$ measured by inertial sensors 
in the vehicle-fixed frame $\sS$; let $\beta$ be the number of such vectors. In addition, consider 
position vectors of a few stationary points in the inertial frame $\sO$ measured by optical (vision 
or lidar) sensors in the vehicle-fixed frame $\sS$. Velocities of the vehicle may be directly measured or can be estimated by linear 
filtering of the optical position vector measurements \cite{Automatica2}. 
Assume that these optical measurements are available for 
$\mpz j$ points at time $t$, whose positions are known in frame $\sO$ as $p_j$, $j\in\cI (t)$, 
where $\cI (t)$ denotes the index set of beacons observed at time $t$. Note that the observed 
stationary beacons or landmarks may vary over time due to the vehicle's motion. These points 
generate ${\mpz j\choose 2}$ unique relative position vectors, which are the vectors connecting 
any two of these landmarks. When two or more position vectors are optically measured, the number of vector measurements that
can be used to estimate attitude is ${\mpz j\choose 2}+\beta$. This number needs to be at least two 
(i.e., ${\mpz j\choose 2}+\beta\geq2$) at an instant, for the attitude to be uniquely determined at 
that instant. In other words, if at least two inertial vectors are measured at all instants (i.e., $\beta\geq2$), 
then beacon position measurements are not required for estimating attitude. However, at least one 
beacon or feature point position measurement is still required to estimate the position of the vehicle.

\subsection{Pose Measurement Model} 
Denote the position of an optical sensor and the unit vector from that sensor to an observed  
beacon in frame $\sS$ as $s^k\in\bR^3$ and $u^k\in\bS^2$, $k=1,\ldots,\mpz k$, respectively. 
Denote the relative position of the $j^{th}$ stationary beacon observed by the $k^{th}$ 
sensor expressed in frame $\sS$ as $q^k_j$. Thus, in the absence of measurement noise
\begin{align}
p_j=R(q^k_j+s^k)+b=Ra_j+b,\; j\in\mathcal I(t),
\label{FrameTrans}
\end{align}
where $a_j=q^k_j+s^k$, are positions of these points expressed in $\sS$. In practice, the $a_j$ 
are obtained from range measurements that have additive noise, which we denote as $a_j^m$. In the case of lidar range measurements, these are given by
\be a_j^m=(q^k_j)^m+s^k=(\varrho_j^k)^m u^k+s^k,\; j\in\mathcal I(t), \label{pim} \ee
where $(\varrho_j^k)^m$ is the measured range to the point by the $k^{th}$ sensor. 
The mean of the vectors $p_j$ and $a_j^m$ are denoted as $\bar{p}$ and $\bar{a}^m$ respectively, 
and satisfy
\begin{align}
\bar{a}^m= R\T (\bar p- b)+ \bar\varsigma, \label{barp}
\end{align}
where $\bar{p}=\frac{1}{\mpz j}\sum\limits^\mpz j_{j=1}p_j$, $\bar{a}^m=\frac{1}{\mpz j}\sum
\limits^\mpz j_{j=1}a_j^m$ and $\bar\varsigma$ is the additive measurement noise 
obtained by averaging the measurement noise vectors for each of the $a_j$. Consider the 
${\mpz j\choose 2}$ relative position vectors from optical measurements, denoted as $d_j=
p_\lambda-p_\ell$ in frame $\sO$ and the corresponding vectors in frame $\sS$ as $l_j=
a_\lambda-a_\ell$, for $\lambda,\ell\in\cI (t)$, $\lambda\ne \ell$. The $\beta$ measured inertial 
vectors are included in the set of $d_j$, and their corresponding measured values expressed 
in frame $\sS$ are included in the set of $l_j$. If the total number of measured (optical and inertial) vectors, ${\mpz j\choose 2}+\beta=2$, then $l_3=l_1\times l_2$ is considered 
a third measured direction in frame $\sS$ with corresponding vector $d_3=d_1\times d_2$ in 
frame $\sO$. Therefore,
\begin{align}
d_j=Rl_j\Rightarrow D=RL,
\end{align}
where $D=[d_1\;\, \cdots\;\, d_n]$, $L=[l_1\;\, \cdots\;\, l_n]\in\bR^{3\times n}$ with $n=3$ if 
${\mpz j\choose 2}+\beta=2$ and $n={\mpz j\choose 2}+\beta$ if ${\mpz j\choose 2}+\beta>2$. 
Note that the matrix $D$ consists of vectors known in frame $\sO$. Denote the measured value 
of matrix $L$ in the presence of measurement noise as $L^m$. Then,
\begin{align}
L^m=R\T D+\mathscr{L},
\label{VecMeasMod}
\end{align}
where $\mathscr{L}\in\bR^{3\times n}$ consists of the additive noise in the vector measurements made
in the body frame $\sS$. 

\subsection{Velocities Measurement Model} 
Denote the angular and translational velocity of the rigid body expressed in frame $\sS$ by $\Omega$ and $\nu$, respectively. Thus, one can write the kinematics 
of the rigid body as
\begin{align}
\dot{\Omega}=R\Omega^\times,\dot{b}=R\nu\Rightarrow\dot{\msg}= \msg\xi^\vee,
\label{Kinematics}
\end{align}
where $\xi= \bbm \Omega\\ \nu\ebm\in\bR^6$ and $\xi^\vee=\bbm \Omega^\times &\; \nu\\ 0 \;\;& 0\ebm$ and $(\cdot)^\times: \bR^3\to\so\subset\bR^{3\times 3}$ is the 
skew-symmetric cross-product operator that gives the vector space isomorphism between 
$\bR^3$ and $\so$. For the general development of the motion estimation scheme, it is assumed that the velocities 
are directly measured. The estimator is then extended to cover the cases where: (i) only angular velocity is 
directly measured; and (ii) none of the velocities are directly measured.

\section{DYNAMIC ESTIMATION OF MOTION FROM PROXIMITY MEASUREMENTS}\label{Sec3}
In order to obtain state estimation schemes from measurements as outlined in Section 
\ref{Sec2} in continuous time, the Lagrange-d'Alembert principle is applied to an action functional of 
a Lagrangian of the state estimate errors, with a dissipation term linear 
in the velocities estimate error. This section presents the estimation scheme obtained 
using this approach. Denote the estimated pose and its kinematics as
\begin{align}
\hat\msg=\bbm \hat R & \;\;\; \hat{b}\\ 0 & \;\;\; 1\ebm\in\SE, \;\;\dot{\hat\msg}=\hat\msg\hat\xi^\vee,
\end{align}
where $\hat\xi$ is rigid body velocities estimate, with $\hat\msg_0$ as the initial pose estimate and the pose estimation 
error as
\begin{align}
\msh=\msg\hat\msg^{-1}=\bbm Q &\;\;\;\; b-Q\hat{b}\\ 0 & \;\;1\ebm=\bbm Q &\;\;\; x\\ 0 &\;\;\;1\ebm\in\SE,
\end{align}
where $Q=R\hat{R}\T$ is the attitude estimation error and $x=b-Q\hat{b}$. Then one obtains, in the case of perfect measurements,
\begin{align}\begin{split}
&\dot\msh = \msh\varphi^\vee,\, \mbox{ where }\, \varphi(\hat\msg,\xi^m,\hat\xi)=\bbm\omega\\\upsilon\ebm= \Ad{\hat\msg}\big(\xi^m- \hat\xi),
\end{split}\label{hdot}
\end{align}
where $\Ad{\mpz{g}}=\bbm \mpz{R}~&~0\\ \mpz{b}^\times\mpz{R}~&~\mpz{R}\ebm$ for $\mpz{g}=\bbm \mpz{R} &\;\; \mpz{b}\\ 0 & \;\;1\ebm$. The attitude and position estimation error dynamics are also in the form
\begin{align}
\dot{Q}=Q\omega^\times,\;\;\dot{x}=Q\upsilon.
\label{Qdot}
\end{align}

\subsection{Lagrangian from Measurement Residuals}
Consider the sum of rotational and translational measurement residuals between the 
measurements and estimated pose as a potential energy-like function. Defining the trace inner 
product on $\mathbb{R}^{n_1\times n_2}$ as
\begin{align}
\lan A_1,A_2\ran :=\tr(A_1 \T A_2),\label{tr_def}
\end{align}
the rotational potential function (Wahba's cost function~\cite{jo:wahba}) is expressed as
\begin{align}
\cU^0_r (\hat{\msg},L^m,D) &= \frac12\lan D -\hat R L^m , (D -\hat R L^m)W\ran,
\label{U0r}
\end{align}
where $W=\diag(w_j)\in\bR^{n\times n}$ is a positive diagonal matrix of weight factors for the 
measured $l_j^m$. Consider the translational potential function
\begin{align}
\cU_t (\hat \msg,\bar a^m,\bar p) &= \frac12 \kappa y\T y= \frac12\kappa \|\bar{p}-\hat{R}\bar{a}^m-\hat{b}\|^2,
\label{U0t} 
\end{align}
where $\bar{p}$ is defined by \eqref{barp}, $y\equiv y(\hat\msg,\bar{a}^m,\bar p)=\bar{p}-
\hat{R}\bar{a}^m-\hat{b}$ and $\kappa$ is a positive scalar. Therefore, the total potential function is 
defined as the sum of the generalization of \eqref{U0r} defined in~\cite{Automatica,ast_acc14}  
for attitude determination on $\SO$, and the translational energy \eqref{U0t} as
\begin{align}
\cU(\hat{\msg},L^m,D,\bar a^m,\bar p)&= \Phi \big( \cU^0_r (\hat{\msg},L^m,D) \big)+\cU_t(\hat{\msg},\bar a^m,\bar p)\nn\\
&=\Phi \big(\frac12\lan D -\hat R L^m , (D -\hat R L^m)W\ran\big)\nn\\
&~~~~~~~+\frac12\kappa \|\bar{p}-\hat{R}\bar{a}^m-\hat{b}\|^2,\label{costU}
\end{align}
where $W$ is positive definite (not necessarily diagonal) which can be selected according to Lemma 3.2 in \cite{Automatica}, and 
$\Phi: [0,\infty)\mapsto[0,\infty)$ is a $\mC^2$ function that satisfies $\Phi(0)=0$ and 
$\Phi'(\mpz x)>0$ for all $\mpz x\in[0,\infty)$. Furthermore, $\Phi'(\cdot)\leq\alpha(\cdot)$ where 
$\alpha(\cdot)$ is a Class-$\mathcal{K}$ function~\cite{khal} and $\Phi'(\cdot)$ denotes the 
derivative of $\Phi(\cdot)$ with respect to its argument. Because of these properties 
of the function $\Phi$, the critical points and their indices coincide for $\cU^0_r$ and 
$\Phi(\cU^0_r)$~\cite{Automatica}. Define the kinetic energy-like function: 
\be
\cT \Big(\varphi(\hat\msg,\xi^m,\hat\xi)\Big)= \frac12 \varphi(\hat\msg,\xi^m,\hat\xi)\T \bJ\varphi(\hat\msg,\xi^m,\hat\xi), 
\label{costT} \ee
where $\bJ\in\bR^{6\times 6}>0$ is an artificial inertia-like kernel matrix. Note that in contrast to 
rigid body inertia matrix, $\bJ$ is not subject to intrinsic physical constraints like the triangle 
inequality, which dictates that the sum of any two eigenvalues of the inertia matrix has to be 
larger than the third. Instead, $\bJ$ is a gain matrix that can be used to tune the estimator. For 
notational convenience, $\varphi(\hat\msg,\xi^m,\hat\xi)$ is denoted as $\varphi$ from 
now on; this quantity is the velocities estimation error in the absence of measurement 
noise. Now define the Lagrangian 
\be
\cL (\hat{\msg},L^m,D,\bar a^m,\bar p,\varphi)= \cT(\varphi) -\cU(\hat{\msg},L^m,D,\bar a^m,\bar p),
\label{contLag}\ee
and the corresponding action functional over an arbitrary time interval $[t_0,T]$ for $T>0$,
\be \cS \big(\cL (\hat{\msg},L^m,D,\bar a^m,\bar p,\varphi)\big)= \int_{t_0}^T \cL 
(\hat{\msg},L^m,D,\bar a^m,\bar p,\varphi) \di t, \,\label{action} \ee
such that $\dot{\hat \msg}= \hat \msg(\hat\xi)^\vee$. A Rayleigh dissipation term linear in the velocities of the form $\bD\varphi$ 
where $\bD\in\bR^{6\times 6}>0$ is used in addition to the Lagrangian \eqref{contLag}, and 
the Lagrange-d'Alembert principle from variational mechanics is applied to obtain the estimator on 
$\Ta\SE$. This yields
\begin{align}
\delta_{\msh,\varphi} \cS \big(\cL (\msh,D,\bar{p},\varphi)\big)=\int_{t_0}^T \eta\T
\bD\varphi \di t,\label{LagdAlem}
\end{align}
which in turn results in the following continuous-time filter.

\subsection{Variational Estimator for Pose and Velocities}
The nonlinear variational estimator obtained by 
applying the Lagrange-d'Alembert principle to the Lagrangian \eqref{contLag} with a dissipation 
term linear in the velocities estimation error, is given by the following statement.

\begin{theorem}\label{filterTHM}
The nonlinear variational  estimator for pose and velocities is given by
\begin{align}
\begin{cases}
\bJ\dot{\varphi}&=\adast{\varphi}\bJ\varphi-Z(\hat{\msg},L^m,D,\bar{a}^m,\bar{p})-\bD\varphi,\vspace{.05in}\\
\hat{\xi}&=\xi^m-\Ad{\hat\msg^{-1}}\varphi,
\vspace{.05in}\\
\dot{\hat{\msg}}&=\hat{\msg}(\hat{\xi})^\vee,\label{ContFil}
\end{cases}
\end{align}
where $\adast{\zeta}=(\ad{\zeta})\T$ with $\ad{\zeta}$ defined by
\begin{align}
\ad{\mpz{\zeta}}=\bbm \mpz w^\times~~ & 0\\ \mpz v^\times\;\; & \mpz{w}^\times\ebm\mbox{ for } \zeta=\bbm \mpz w\\ \mpz v\ebm\in\bR^6,\label{ad_def}
\end{align}
and $Z(\hat{\msg},L^m,D,\bar{a}^m,\bar{p})$ is defined by
\begin{align}
\begin{split}
Z(\hat{\msg},L^m,D,&\bar{a}^m,\bar{p})=\label{Z}\\
&\bbm \Phi'\Big(\cU^0_r (\hat{\msg},L^m,D)\Big)S_\Gamma(\hat{R})+\kappa\bar{p}^\times y\\ \kappa y\ebm,\end{split}
\end{align}
where $\cU^0_r (\hat{\msg},L^m,D)$ is defined as \eqref{U0r}, $y\equiv y(\hat{\msg},\bar{a}^m,\bar{p})=\bar{p}-\hat{R}\bar{a}^m-\hat{b}$ and
\begin{align} S_\Gamma(\hat{R})=\mrm{vex}\big(DW(L^m)\T\hat{R}\T-\hat{R}L^mWD\T\big), \label{SLdef} \end{align}
where $\mrm{vex}(\cdot): \so\to\bR^3$ is the inverse of the $(\cdot)^\times$ map.
\end{theorem}
The proof is presented in \cite{Automatica2,Gaurav_ASR}. In the proposed approach, the time 
evolution of $(\hat \msg,\hat\xi)$ has the form of the dynamics of a rigid body with Rayleigh 
dissipation. This results in an estimator for the motion states $(\msg,\xi)$ that 
dissipates the ``energy" content in the estimation errors $(\mathsf{h},\varphi)= (\msg \hat\msg^{-1}, 
\Ad{\hat\msg}(\xi- \hat\xi))$ to provide guaranteed asymptotic stability in the case of perfect 
measurements~\cite{Automatica}.

Explicit expressions for the vector of velocities $\xi^m$ can be obtained for two common cases 
when these velocities are not directly measured. These cases are presented here.

\subsection{Variational Estimator Implemented without Direct Velocity Measurements}\label{ContButterworth}
The velocity measurements in \eqref{ContFil} can be replaced by filtered velocity estimates 
obtained by linear filtering of optical and inertial measurements using, e.g., a second-order 
Butterworth filter. This is both useful and necessary when velocities are not directly measured.

\subsubsection{Angular velocity is measured using rate gyros}
For the case that angular velocities are measured by rate gyros besides the 
$\mpz{j}$ feature point position measurements, the linear velocities of the rigid body can 
be calculated using each single position measurement by rewriting \eqref{pjdot} as
\begin{align}
\nu^f=(a_j^f)^\times\Omega^f-v_j^f,
\end{align}
for the $j^{th}$ point. Averaging the values of $\nu$ derived from all feature points gives a more reliable result. Therefore, the rigid body's filtered velocities are expressed in this case as
\begin{align}
\xi^f=\bbm \Omega^f\\ \frac{1}{\mpz{j}}\sum\limits_{j=1}^{\mpz{j}}(a_j^f)^\times\Omega^f-v_j^f \ebm.
\end{align}

\subsubsection{Translational and angular velocity measurements are not available}
In this case, rigid body velocities can be calculated in terms of the measurements. In order to do 
so, one can differentiate \eqref{FrameTrans} as follows
\begin{align}
&\dot{p}_j=R\Omega^\times a_j+R\dot{a}_j+\dot{b}=R\big(\Omega^\times a_j+\dot{a}_j+\nu\big)=0\nn\\
\Rightarrow&\dot{a}_j-a_j^\times\Omega+\nu=0\nn\\
\Rightarrow&v_j=\dot{a}_j=[a_j^\times\; -I]\xi=G(a_j)\xi,
\label{pjdot}
\end{align}
where $G(a_j)=[a_j^\times\; -I]$ has full row rank. From vision-based or Doppler lidar sensors, one can also 
measure the velocities of the observed points in frame $\sS$, denoted $v_i^m$. Here, velocity 
measurements as would be obtained from vision-based sensors is considered. 
The measurement model for the velocity is of the form
\be v_j^m=G(a_j)\xi+\vartheta_j,\ee
where $\vartheta_j\in\bR^3$ is the additive error in velocity measurement $v_j^m$. Instantaneous 
angular and translational velocity determination from such measurements is treated 
in~\cite{ast_acc14}. 
Note that $v_j=\dot{a}_j$, for $j\in\mathcal I(t)$. As this kinematics indicates, the relative 
velocities of at least three beacons are needed to determine the vehicle's translational and 
angular velocities uniquely at each instant. The rigid body velocities are obtained 
using the pseudo-inverse of $\mathds{G}(A^f)$:
\begin{align}
\mathds{G}(A^f)\xi^f&=\mathds{V}(V^f)\Rightarrow\xi^f=\mathds{G}^\ddag(A^f)\mathds{V}(V^f),\label{ximMore2}\\
\mbox{where }\;\mathds{G}(A^f)&=\bbm G(a^f_1)\\\vdots\\G(a^f_\mpz j)\ebm\;\mbox{and }\;\mathds{V}(V^f)=\bbm v^f_1\\\vdots\\v^f_\mpz j \ebm,\label{GVDef}
\end{align}
for $1,...,\mpz j\in\mathcal I(t)$. 
When at least three beacons are measured, $\mathds{G}(A^f)$ is a full column rank matrix, 
and $\mathds{G}^\ddag(A^f)= \Big( \mathds{G}\T(A^f)\mathds{G} (A^f)\Big)^{-1} 
\mathds{G}\T (A^f)$ gives its pseudo-inverse.

\section{DISCRETIZATION FOR COMPUTER IMPLEMENTATION}\label{Sec5}
For onboard computer implementation, the variational estimation scheme outlined above has to 
be discretized. Since the estimation scheme proposed here is obtained from 
a variational principle of mechanics, it can be discretized by applying the discrete 
Lagrange-d'Alembert principle~\cite{marswest}. Consider an interval of time $[t_0, T]\in\bR^+$ 
separated into $N$ equal-length subintervals $[t_i,t_{i+1}]$ for $i=0,1,\ldots,N$, with $t_N=T$ 
and $t_{i+1}-t_i=\Delta t$ is the time step size. Let $(\hat \msg_i,\hat\xi_i)\in\SE\times\bR^6$ 
denote the discrete state estimate at time $t_i$, such that $(\hat \msg_i,\hat\xi_i)\approx 
(\hat \msg(t_i),\hat\xi(t_i))$ where $(\hat \msg(t),\hat\xi(t))$ is the exact solution of the 
continuous-time estimator at time $t\in [t_0, T]$. Let the values of the discrete-time measurements  
$\xi^m$, $\bar a^m$ and $L^m$ at time $t_i$ be denoted as $\xi^m_i$, $\bar a^m_i$ and $L^m_i$, 
respectively. Further, denote the corresponding values for the latter two quantities in inertial frame 
at time $t_i$ by $\bar p_i$ and $D_i$, respectively. The discrete-time filter is then presented in the form of a Lie group variational integrator (LGVI) in the following statement.

\begin{theorem} \label{discfilter}
A first-order discretization of the estimator proposed in Theorem \ref{filterTHM} is given by
\begin{align}
(J\omega_i)^\times&=\frac{1}{\Delta t}(F_i\mathcal{J}-\mathcal{J}F_i\T),\label{LGVI_F}\\
(M+\Delta t\bD_t)\upsilon_{i+1}&=F_i\T M\upsilon_i\label{LGVI_upsilon}\\
&~~~~~~~~+\Delta t \kappa (\hat{b}_{i+1}+\hat{R}_{i+1}\bar{a}^m_{i+1}-\bar{p}_{i+1}),\nn\\
(J+\Delta t\bD_r)\omega_{i+1}&=F_i\T J\omega_i+\Delta t M\upsilon_{i+1}\times\upsilon_{i+1}\nn\\
+\Delta t&\kappa \bar{p}_{i+1}^\times (\hat{b}_{i+1}+\hat{R}_{i+1}\bar{a}^m_{i+1})\label{LGVI_omega}\\
-\Delta t&\Phi' \big( \cU^0_r (\hat{\msg}_{i+1},L_{i+1}^m,D_{i+1}) \big)S_{\Gamma_{i+1}}(\hat{R}_{i+1}),\nn\\
\hat\xi_i&=\xi^m_i-\Ad{\hat\msg_i^{-1}}\varphi_i,\label{LGVI_xihat}\\
\hat\msg_{i+1}&=\hat\msg_i\exp(\Delta t\hat\xi_i^\vee),\label{LGVI_ghat}
\end{align}
where $F_i\in\SO$, $\big(\hat\msg(t_0),\hat\xi(t_0)\big)=(\hat\msg_0,\hat\xi_0)$, $\mathcal{J}$ is defined in terms of positive
matrix $J$ by $\mathcal{J}=\frac12\tr[J]I-J$, $M$ is a positive definite matrix,
$\varphi_i=[\omega_i\T\;\upsilon_i\T]\T$, and $S_{\Gamma_i}(\hat R_i)$ is the value of  
$S_\Gamma(\hat R)$ at time $t_i$, with $S_\Gamma(\hat R)$ defined by \eqref{SLdef}.
\end{theorem}
The proof is presented in \cite{Automatica2}.

\section{NUMERICAL SIMULATIONS}\label{Sec6}
\begin{figure}
\begin{center}
\includegraphics[height=2.6in]{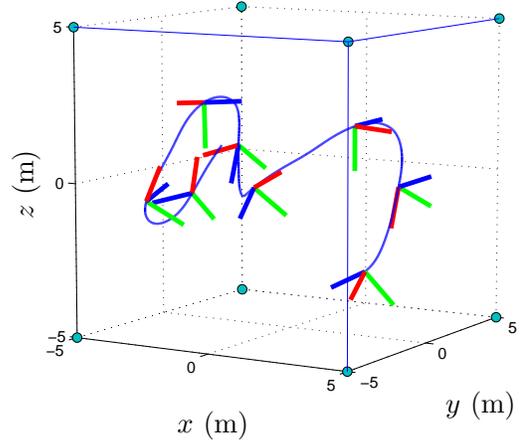}
\caption{Position and attitude trajectory of the simulated vehicle.}  
\label{Fig2}                                 
\end{center}                                 
\end{figure}

\begin{figure}
\begin{center}
\includegraphics[height=2.4in]{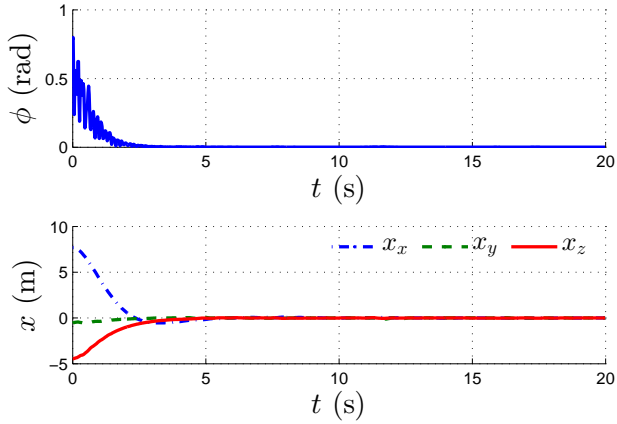}
\caption{Principal angle of the attitude and position estimation error for CASE 1.}  
\label{Fig3}                                 
\end{center}                                 
\end{figure}

\begin{figure}
\begin{center}
\includegraphics[height=2.4in]{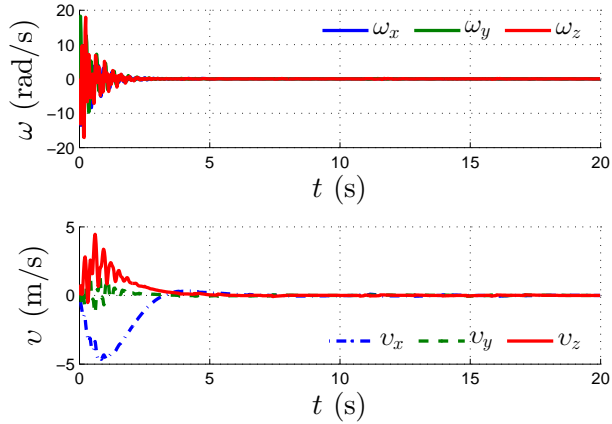}
\caption{Angular and translational velocity estimation error for CASE 1.}  
\label{Fig4}                                 
\end{center}                                 
\end{figure}

This section presents numerical simulation results for the discrete-time estimator obtained in 
Section \ref{Sec5}. In order to numerically simulate this estimator, simulated true states of an 
aerial vehicle flying in a cubical room are produced using the kinematics and dynamics equations of a 
rigid body. The vehicle mass and moment of inertia are taken to be $m_v=420$ g and $J_v=
[51.2\;\;60.2\;\;59.6]\T$ g.m$^2$, respectively. The resultant external forces and torques applied 
on the vehicle are $\phi_v(t)=10^{-3}[10\cos(0.1t)\;\;2\sin(0.2t)\;\;-2\sin(0.5t)]\T$ N and 
$\tau_v(t)=10^{-6}\phi_v(t)$ N.m, respectively. The room is assumed to be a cube of size 
10m$\times$10m$\times$10m with the inertial frame origin at the geometric center. The initial 
attitude and position of the vehicle are:
\begin{align}
R_0&=\expm_{\SO}\bigg(\Big(\frac{\pi}{4}\times[\frac{3}{7}\ -\frac{6}{7}\ \frac{2}{7}]\T\Big)^\times
\bigg),\nn\\
\mbox{and } b_0&=[2.5\ 0.5\ -3]\T\mbox{ m}.
\end{align}
This vehicle's initial angular and translational velocity respectively, are:
\begin{align}
\begin{split}
\Omega_0&=[0.2\;\; -0.05\;\; 0.1]\T\mbox{ rad/s},\\
\mbox{and } \nu_0&=[-0.05\;\;0.15\;\;0.03]\T\mbox{ m/s}.
\end{split}
\end{align}
The vehicle dynamics is simulated over a time interval of $T=150 \mbox{ s}$, with a time 
stepsize of $\Delta t=0.02 \mbox{ s}$. The trajectory of the vehicle over this time interval is 
depicted in Fig.~\ref{Fig2}.
The following two inertial directions, corresponding to nadir and Earth's magnetic field 
direction, are measured by the inertial sensors on the vehicle:
\begin{align}
d_1=[0\;\;0\;\;-1]\T,\;\;d_2=[0.1\;\;0.975\;\;-0.2]\T.
\end{align}
For optical measurements, eight beacons are located at the eight vertices of the cube, labeled 
1 to 8. The positions of these beacons are known in the inertial frame and their index (label) 
and relative positions are measured by optical sensors onboard the vehicle whenever the 
beacons come into the field of view of the sensors. Three identical cameras (optical sensors) 
and inertial sensors are assumed to be installed on the vehicle. The cameras are fixed to known 
positions on the vehicle, on a hypothetical horizontal plane passing through the vehicle, 120$^\circ$ 
apart from each other, as shown in Fig.~\ref{Frames}. All the camera readings contain random 
zero mean signals whose probability distributions are normalized bump functions with width of 
$0.001$m. The following are selected for the positive definite estimator gain matrices:
\begin{align}
J&=\diag\big([0.9\;\;0.6\;\;0.3]\big), \nn \\
M&=\diag\big([0.0608\;\;0.0486\;\;0.0365]\big), \\
\bD_r&=\diag\big([2.7\;2.2\;1.5]\big),\bD_t=\diag\big([0.1\;\;0.12\;\;0.14]\big). \nn
\end{align}
$\Phi(\cdot)$ could be any $C^2$ function with the properties described in Section \ref{Sec3}, but is 
selected to be $\Phi(x)=x$ here. The initial state estimates have the following values:
\begin{align}
\begin{split}
\hat\msg_0&=I,\;\;\;\hat\Omega_0=[0.1\;\;0.45\;\;0.05]\T\mbox{ rad/s},\\
\mbox{ and }\hat\nu_0&=[2.05\;\;0.64\;\;1.29]\T\mbox{ m/s}.
\end{split}
\end{align}
A conic field of view (FOV) of 2$\times$40$^\circ$ is assumed for the cameras, which guarantees 
at least three beacons observed are common between successive measurements. The vehicle's 
velocity vector is calculated from \eqref{ximMore2}. The discrete-time estimator 
\eqref{LGVI_F}-\eqref{LGVI_ghat} is simulated over a time of $T=20$ s with time stepsize 
$\Delta t=0.02$ s. At each measurement instant, \eqref{LGVI_F} is solved using Newton-Raphson 
iterations to find an approximation for $F_i$. The remaining equations (all explicit) are 
solved consecutively to generate the estimated states. The principal angle of the attitude 
estimation error and the position estimate error are plotted in Fig.~\ref{Fig3}. The angular 
and translational components of the vehicle's velocity estimate errors are also depicted in 
Fig.~\ref{Fig4}.

\section{CONCLUSION}\label{Sec7}
This article proposes an estimator for rigid body pose and velocities, using optical and inertial 
measurements by sensors onboard the rigid body. The sensors are assumed to provide 
measurements in continuous-time or at a sufficiently high frequency, with bounded noise. An 
artificial kinetic energy quadratic in rigid body velocity estimate errors is 
defined, as well as two fictitious potential energies: (1) a generalized Wahba's cost function for 
attitude estimation error in the form of a Morse function, and (2) a quadratic function of  
the vehicle's position estimate error. Applying the Lagrange-d'Alembert principle on a Lagrangian 
consisting of these energy-like terms and a dissipation term linear in velocities estimation error, an 
estimator is designed on the Lie group of rigid body motions. A discrete-time counterpart of this estimator is also presented. 
In the presence of measurement noise, numerical simulations show that state estimates converge 
to a bounded neighborhood of the true states.
\pagestyle{empty}

\addtolength{\textheight}{-2cm}   

\end{document}